\documentclass[a4paper,nopreprintline]{elsarticle}
\usepackage{hyperref}
\usepackage{amsmath,amssymb,amsthm}
\usepackage{enumerate,xcolor}
\numberwithin{equation}{section}
\usepackage{shortvrb, makeidx}

\usepackage[pagewise]{lineno}
\AtBeginDocument{%
  \abovedisplayskip=0.5\abovedisplayskip
  \belowdisplayskip=0.5\belowdisplayskip
  \abovedisplayshortskip=0.5\abovedisplayshortskip
  \belowdisplayshortskip=0.5\belowdisplayshortskip
}

\newtheorem{definition}{\hspace{2em}Definition}[section]
\newtheorem{theorem}[definition]{\hspace{2em}Theorem}

\newtheorem{lemma}[definition]{\hspace{2em}Lemma}
\newtheorem{proposition}[definition]{\hspace{2em}Proposition}

\newtheorem{remark}{\hspace{2em}Remark}[section]

\newcommand{\rn}{\int_{\mathbb{R}^N}}

\everymath{\displaystyle}
\begin{document}

\title{Ground state and multiple solutions for a class of modified autonomous fourth-order elliptic equations under Berestycki--Lions type conditions}
\author[LFY]{Lifeng Yin}
\ead{yin136823@163.com}
\author[FW]{Fan ~Wang\corref{cor1}}
\ead{wangf767@swjtu.edu.cn}
\cortext[cor1]{Corresponding author.}
\address[LFY]{School of Mathematical Sciences, Sichuan Normal University, Chengdu 610066, China.}
\address[FW]{Department of Mathematics, Southwest Jiaotong University, Chengdu 611756, China.}
\date{}

\begin{abstract}
This paper investigates the following modified fourth-order elliptic equation
    \[
    \begin{aligned}
        \left\{
        \begin{array}{ll}
            \Delta^2 u - \Delta u + u - \frac{1}{2}u\Delta(u^2) = f(u), & \text{in } \mathbb{R}^N, \\
            u \in H^2(\mathbb{R}^N),
        \end{array}
        \right.
    \end{aligned}
    \]
    where $N=5$ or $6$, and the nonlinearity $f:\mathbb{R}\rightarrow\mathbb{R}$ satisfies Berestycki--Lions type conditions.
    By using Jeanjean's augmented functional technique and constructing a Poho\v{z}aev--Palais--Smale sequence, we prove that the problem admits a ground state solution. Moreover, if \(f\) is odd, we prove, via a symmetric mountain-pass lemma, that the problem admits infinitely many radially symmetric solutions.
    \end{abstract}

\begin{keyword}
Biharmonic operator; Poho\v{z}aev identity; Multidimensional odd paths; Ground state and multiple solutions
\end{keyword}

\maketitle

%The title of your section 1
\section{Introduction}
In this paper, we consider the following fourth-order elliptic equation
\begin{align}\label{eq1}
\left\{
\begin{array}{ll}
 \Delta^2 u-\Delta u+u-\frac{1}{2}u\Delta(u^2)=f(u), & \text{in } \mathbb{R}^{N}, \\
 u\in H^2(\mathbb{R}^N),
 \end{array}
 \right.
\end{align}
where $N=5$ or $6$ and $f:\mathbb{R}\to\mathbb{R}$ is continuous and satisfies the conditions specified below.

Equation \eqref{eq1} is a special case of the following modified quasilinear fourth-order elliptic problem:
\begin{align}\label{eq2}
    \Delta^2 u -\lambda  \Delta u + V(x)u - \kappa u \Delta(u^2) = g(x, u), \quad \text{in } \mathbb{R}^{N},
\end{align}
with $\kappa \geq 0$ and $\lambda\in\mathbb{R}$. When $\mathbb{R}^N$ in \eqref{eq2} is replaced by a bounded domain $\Omega$ and $\kappa = 0$, the equation models traveling waves in suspension bridges, as studied by Lazer and McKenna \cite{MR1084570} and McKenna and Walter \cite{MR1050908}.

Another related topic is the normalized biharmonic Schr\"odinger equation
\begin{align}\label{eq3}
  \left\{
  \begin{array}{ll}
    \Delta^2 u + \lambda u = h(x,u), & x\in\mathbb R^N, \\
    \int_{\mathbb{R}^N} u^2 \, dx = c,
  \end{array}
  \right.
\end{align}
where $\lambda$ is the Lagrange multiplier associated with the prescribed mass.  Liu et al. \cite{Liu2025} considered a spatially dependent nonlinearity
\[h(x,u) = a(\varepsilon x)|u|^{q-2}u + |u|^{p-2}u\text{,}\]
where $a:\mathbb{R}^N\rightarrow[0,\infty)$ is continuous and satisfies suitable conditions. Assuming $2 < q < 2 + \frac{8}{N} < p < 2_* := \frac{2N}{N-4}$ and $N \geq 5$, they established the existence of multiple normalized solutions of \eqref{eq3} via variational methods. Zhang et al. \cite{Zhang2025} extended these results to nonlinearities $h(x,u)=h(u)$ satisfying:
\begin{itemize}
\item[(H1)] $h: \mathbb{R}\rightarrow\mathbb{R}$ is continuous and odd;
    \item[(H2)]  There exist constants $2 + \frac{8}{N} < \alpha \leq \beta < 2_*$ such that
    \[
    0 < \alpha H(s) \leq h(s)s \leq \beta H(s), \quad \forall s \neq 0,
    \]
    where $H(s) = \int_0^s h(t)  dt$;
    \item[(H3)]  The function defined by
$
\tilde{H}(s) := \frac{1}{2}h(s)s - H(s)
$
is of class \(C^1\) and satisfies
\[
\tilde{H}'(s)s \geq \alpha \tilde{H}(s), \quad \forall s \in \mathbb{R},
\]
where \(\alpha\) is given in (H2).
\end{itemize}
Their novel constrained variational methods yielded normalized ground state solutions of \eqref{eq3} under general mass-supercritical assumptions.

Mederski and Siemianowski \cite{Mederski2023} considered the autonomous case
\begin{align*}
  \Delta^2 u = h(u), \quad x \in \mathbb{R}^N  ,
\end{align*}
 where $N \geq 5$. Using the Poho\v{z}aev manifold method, they proved the existence of ground states for this equation when $h$ has general subcritical growth at infinity, without assuming the Ambrosetti--Rabinowitz type condition. Further contributions can be found in \cite{Alarcon2025,Fernandez2022,Lu2024,Yang2025} and the references therein.

In the non-autonomous case ($V \not\equiv \mathrm{constant}$), Chen, Liu, and Wu \cite{MR3276713} proved the existence of nontrivial and high-energy solutions for \eqref{eq2} under coercive potentials and for $3$-superlinear nonlinearities $g$ satisfying:
 \[\lim_{|t|\rightarrow+\infty}\frac{G(x,t)}{t^4}=+\infty,\qquad G(x,t):=\int_0^t g(x,s)\,ds,\]
  together with other mild assumptions. These results were extended by Cheng and Tang \cite{MR3591225}. Using sign-changing Nehari manifolds, the authors of \cite{MR3348950} obtained least-energy nodal solutions of \eqref{eq2} in the critical-growth case. Liu and Zhao \cite{MR3945610} used Morse theory to obtain nontrivial solutions of \eqref{eq2}, while Yin and Jiang \cite{MR4291515} studied indefinite potentials with nonlinearities $g(x,u)=g(u)$ satisfying:
\begin{itemize}
    \item[$(g)$] There exists $\mu > 3$ such that $g(t)t \geq \mu \int_0^t g(s)  ds > 0$  for all $t \neq 0$.
\end{itemize}

Motivated by these works, we establish ground state and multiplicity results for the autonomous problem \eqref{eq1}. To the best of our knowledge, no prior results exist for \eqref{eq1} under Berestycki--Lions type conditions \cite{MR695535}.

 To state our main theorems, we assume that $f$ satisfies the following conditions:
\begin{enumerate}[\indent($f_1$)]
\item  $f\in C(\mathbb R)$, $\lim_{t\rightarrow0}\frac{f(t)}{t}=0$;
\item $\lim_{|t|\rightarrow\infty}\frac{f(t)}{|t|^{2_*-1}}=0$, where $2_*:=\frac{2N}{N-4}$;
\item there exists $\xi>0$ such that $\widetilde{F}(\xi):=F(\xi)-\frac{1}{2}\xi^2>0$, where $F(t):=\int_0^t f(s)\,ds$ for $t\in\mathbb R$.
\end{enumerate}

Our main results are:
\begin{theorem}\label{thm1}
 Assume that $\left(f_{1}\right)$--$\left(f_{3}\right)$ hold. Then equation \eqref{eq1} admits at least one ground-state solution.
\end{theorem}

\begin{remark}
  The hypotheses $(f_1)$--$(f_3)$ are of Berestycki--Lions type and, in particular, do not impose the Ambrosetti--Rabinowitz condition.
\end{remark}

\begin{theorem}\label{th2}
Assume that $\left(f_{1}\right)$--$\left(f_{3}\right)$ hold and that $f$ is odd. Then equation \eqref{eq1} admits infinitely many radially symmetric solutions.
\end{theorem}

The proofs of Theorems \ref{thm1} and \ref{th2} employ variational methods. For Theorem \ref{thm1}, the main difficulty in studying \eqref{eq1} when $f$ satisfies only the Berestycki--Lions conditions is that it is not known whether Palais--Smale sequences (or even Cerami) sequences are bounded or not. Typically, such boundedness requires the Ambrosetti--Rabinowitz condition:
\begin{enumerate}
  \item[$(\ast)$] $\exists \mu>4$ such that $f(t)t\geq \mu F(t)\geq 0$ for $t\in\mathbb{R}$.
\end{enumerate}
To circumvent this, we adapt Jeanjean's technique \cite{MR1430506} through:
\begin{itemize}
    \item defining the scaling map $\varphi: E\rightarrow H^2(\mathbb{R}^N)$, $\varphi(s,v)(x) = v(x/e^s)$;
    \item constructing the augmented functional $\widetilde{I}=I\circ \varphi$;
    \item generating a Poho\v{z}aev--Palais--Smale sequence $\{u_n\} \subset H^2(\mathbb{R}^N)$ at level $c$ (see \eqref{pps}).
\end{itemize}
This framework yields a bounded Palais--Smale sequence $\{u_n\}$ for $I$.

For Theorem \ref{th2}, we utilize the comparison functional $J$ from \eqref{cf}. A key ingredient is the construction of
\[ \tau: \pi_{k-1} \to H_{0,r}^2(B_R) \]
satisfying the following properties:
\begin{align*}
0\not\in \tau(\pi_{k-1}), \int_{\mathbb{R}^N}\widetilde{F}(\tau(\ell))dx\geq 1 \quad \forall\ell\in \pi_{k-1}.
\end{align*}
This construction yields two essential results:
\begin{enumerate}[$(i)$]
    \item Both $I$ and $J$ satisfy the symmetric mountain-pass geometry;
    \item The minimax classes in \eqref{eq:3.1}--\eqref{eq:3.2} are nonempty, and the associated levels
$b_n$ and $c_n$ are well defined.  The criticality of $c_n$ is addressed in
Proposition~\ref{p2}, while that of $b_n$ is proved later by means of the augmented functional.

%
%    The critical values $b_n$ and $c_n$ defined in \eqref{eq:3.1} are well-defined (Proposition \ref{p2}).
\end{enumerate}
The subsequent boundedness arguments follow the same strategy as in the proof of Theorem \ref{thm1}.

This paper is organized as follows: Section 2 presents preliminaries, Section 3 proves Theorem \ref{thm1}, and Section 4 establishes Theorem \ref{th2} via minimax methods.

\section{Preliminaries}
Throughout this paper, we use the following notation:
\begin{itemize}
 \item $L^p(\mathbb{R}^N)$: Lebesgue space with norm $\|u\|_p = \left(\int_{\mathbb{R}^N}|u|^pdx\right)^{1/p}$.
 \item $H^{2}(\mathbb{R}^N) = \{u \in L^2(\mathbb{R}^N) : \Delta u, \nabla u \in L^2(\mathbb{R}^N)\}$ endowed with the norm $\|u\| = \left( \int_{\mathbb{R}^N} (|\Delta u|^2 + |\nabla u|^2 + u^2) dx \right)^{1/2}$.
 \item $H_r^2(\mathbb{R}^N)$: subspace of radially symmetric functions in $H^2(\mathbb{R}^N)$.
 \item 
 \(
 H^2_{0,r}(B_R)
 :=H^2_0(B_R)\cap\{u:u\text{ is radial}\}.
\)
 \item $D^{2,2}(\mathbb{R}^N)$: completion of $C_0^\infty(\mathbb{R}^N)$ under $\|u\|_{D^{2,2}} = \|\Delta u\|_2$.
 \item $E = \mathbb{R} \times H^2(\mathbb{R}^N)$ with norm $\|(s,u)\|_E = (|s|^2 + \|u\|^2)^{1/2}$.
 \item $E'$: dual space of $E$.
% \item $B_r(y) = \{x \in \mathbb{R}^N : |x-y| \leq r\}$, $B_r = B_r(0)$.
 \item $u_t(x) = u(x/t)$ for any $t > 0$.
 \item $C > 0$: generic constant that may vary between estimates.
\end{itemize}
Since $N=5$ or $6$, the Sobolev embeddings $H^2(\mathbb R^N)\hookrightarrow L^6(\mathbb R^N)$ and $H^2(\mathbb R^N)\hookrightarrow W^{1,3}(\mathbb R^N)$ hold. Hence, following the proof of \cite[Lemma 2.2]{MR3276713}, we have the following estimate for all $u\in H^{2}\left(\mathbb{R}^{N}\right)$:
\begin{align}\label{1}
\int_{\mathbb{R}^N} u^2 |\nabla u|^2 dx \leq \|u\|_6^2 \|\nabla u\|_3^2 \leq \|u\|_6^2 \|u\|_{W^{1,3}}^2 \leq C \|u\|^4.
\end{align}
The energy functional $I: H^2(\mathbb{R}^N) \to \mathbb{R}$ defined by
\begin{align*}
I(u) = \frac{1}{2} \int_{\mathbb{R}^{N}} (|\Delta u|^2 + |\nabla u|^{2} + u^{2}) dx + \frac{1}{2} \int_{\mathbb{R}^N} u^2 |\nabla u|^2 dx - \int_{\mathbb{R}^{N}} F(u) dx
\end{align*}
is well-defined and of class $C^1$ under assumptions $(f_1)$--$(f_2)$ and \eqref{1}. Its Fr\'{e}chet derivative is given by
\begin{align*}
\langle I'(u), v \rangle = &\int_{\mathbb{R}^{N}} \left( \Delta u \Delta v + \nabla u \cdot \nabla v + u v + uv|\nabla u|^2 + u^2 \nabla u \cdot \nabla v \right) dx\notag\\
&\quad - \int_{\mathbb{R}^{N}} f(u)v dx
\end{align*}
for $v \in H^{2}(\mathbb{R}^{N})$.

\begin{lemma}\label{lem:mountain_pass_geometry}
Under $(f_1)$--$(f_3)$, $I$ satisfies the mountain-pass geometry:
\begin{enumerate}[(i)]
  \item There exist $\rho, b > 0$, with $\rho$ sufficiently small, such that $\inf_{\|u\| = \rho} I(u) \geq b$,\label{mp1}
  \item There exists $ u_0 \in H^2(\mathbb{R}^N)$ with $\|u_0\| > \rho$ and $I(u_0) < 0$.\label{mp2}
\end{enumerate}
\end{lemma}

\begin{proof}
(\ref{mp1}) From $(f_1)$--$(f_2)$, for any $\varepsilon > 0$ there exists $C_\varepsilon > 0$ such that
\begin{equation}\label{f1}
|f(u)| \leq \varepsilon |u| + C_\varepsilon |u|^{2_*-1}.
\end{equation}
Integrating \eqref{f1} gives $|F(u)|\le \frac{\varepsilon}{2}u^2+C_\varepsilon |u|^{2_*}$. Using this inequality with $\varepsilon=\frac{1}{2}$ and  the Sobolev inequality, there exist $b,\rho>0$ such that
\begin{align*}
I(u) &\geq \frac{1}{2} \|u\|^2 - \int_{\mathbb{R}^N} F(u) dx \\
&\geq \frac{1}{2} \|u\|^2 - \frac{1}{4} \|u\|_2^2 - C \|u\|_{2_*}^{2_*} \\
&\geq \frac{1}{4} \|u\|^2 - C_1 \|u\|^{2_*}\\
&=\|u\|^2\left(\frac{1}{4}- C_1 \rho^{2_*-2}\right)=b > 0
\end{align*}
for $\|u\|=\rho$ sufficiently small.

(\ref{mp2}) Following the cut-off construction in \cite{MR695535}, let \(R>1\) be sufficiently large. Define $\eta$ by:
\[
\eta(s)=
\begin{cases}
1, & s\leq 0,\\
1-3s^2+2s^3, & 0<s<1,\\
0, & s\geq 1.
\end{cases}
\]
Set
\begin{equation}\label{eq:test_function}
v(x):=\xi \eta(|x|-R).
\end{equation}
Then \(v\in H^2(\mathbb R^N)\), \(v=\xi\) in \(B_R\),
\(v=0\) in \(\mathbb R^N\setminus B_{R+1}\), and
\(0\leq v\leq \xi\).
By $(f_3)$ and \eqref{eq:test_function}, one has
\begin{equation}\label{vR}
\begin{split}
\int_{\mathbb{R}^N} \widetilde{F}(v) dx &\geq \widetilde{F}(\xi) |B_R| - |B_{R+1} \setminus B_R| \max_{t \in [0,\xi]} |\widetilde{F}(t)| \\
&\geq C R^N - C' R^{N-1} > 0
\end{split}
\end{equation}
for $R$ sufficiently large. By \eqref{vR}, we have
\begin{align*}
\lim_{t \to \infty} \frac{I(v_t)}{t^N} &= \lim_{t \to \infty} \frac{1}{t^N} \left[ \frac{1}{2} \left( t^{N-4} \|\Delta v\|_2^2 + t^{N-2} \|\nabla v\|_2^2 + t^N \|v\|_2^2 \right) \right] \notag\\
&\quad+ \lim_{t \to \infty} \frac{1}{t^N} \left[ \frac{t^{N-2}}{2} \int v^2 |\nabla v|^2 dx - t^N \int F(v) dx \right] \notag\\
&= - \int_{\mathbb{R}^N} \widetilde{F}(v)  dx < 0.
\end{align*}
Thus $I(v_t) \to -\infty$ as $t \to \infty$. Choose $t_0$ sufficiently large so that
\[
 I(v_{t_0})<0,
 \qquad
 \|v_{t_0}\|>\rho,
\]
and set $u_0:=v_{t_0}$.
\end{proof}
Define
\[ Q(u)=\int_{\mathbb R^N}u^2|\nabla u|^2\,dx. \]
By \cite[Lemma 4.1]{MR3945610} and $u_n\rightharpoonup u$ in $H^2(\mathbb{R}^N)$, we have
\begin{align}
  \label{qslinear}
Q(u) \leq \liminf_{n\rightarrow+\infty} Q(u_n), \quad \langle Q'(u_n), \phi \rangle \to \langle Q'(u), \phi \rangle,
\end{align}
for all $\phi\in H^2(\mathbb{R}^N)$.
\begin{lemma}\label{lem:nonvanishing}
Let $\{u_n\} \subset H^2(\mathbb{R}^N)$ be a bounded Palais--Smale sequence for $I$ at a positive level $c>0$. Then the set of critical points
$\mathcal{N} = \{u \in H^2(\mathbb{R}^N) \setminus \{0\} : I'(u) = 0\}$ is nonempty.
\end{lemma}

\begin{proof}
Since $\{u_n\}$ is bounded in $H^2(\mathbb R^N)$, after passing to a subsequence,
$u_n\rightharpoonup u$ in $H^2(\mathbb R^N)$ and $u_n\to u$ in
$L^p_{\rm loc}(\mathbb R^N)$ for every $p\in[2,2_*)$. Fix $r>0$ and suppose that
\begin{align}\label{s1}
     \lim_{n\to\infty}\sup_{y\in\mathbb R^N}\int_{B_r(y)}|u_n|^2dx=0 .
\end{align}
Then Lions' vanishing lemma yields
\begin{align}\label{f6}
     u_n\to0\quad\hbox{in }L^p(\mathbb R^N)\quad\hbox{for all }p\in(2,2_*).
\end{align}
From $(f_1)$--$(f_2)$, for any $\varepsilon>0$, there exists $C_\varepsilon>0$ such that
\begin{align}\label{f2}
  |f(s)|\leq \varepsilon\bigl(|s|+|s|^{2_*-1}\bigr)+C_\varepsilon |s|^{q-1},\quad \text{for~~}q\in(2,2_*).
\end{align}
Using \eqref{f6}, the boundedness of $\{u_n\}$ in $L^{2_*}(\mathbb R^N)$ and the arbitrariness of $\varepsilon$, we obtain
\[
 \int_{\mathbb R^N} f(u_n)u_n\,dx=o(1),\qquad
 \int_{\mathbb R^N} F(u_n)\,dx=o(1).
\]
Therefore
\[
\begin{split}
 I(u_n)-\frac12\langle I'(u_n),u_n\rangle
 &=\int_{\mathbb R^N}\left(\frac12 f(u_n)u_n-F(u_n)\right)dx
   -\frac12\int_{\mathbb R^N}u_n^2|\nabla u_n|^2dx  \\
 &\leq o(1).
\end{split}
\]
On the other hand, since $I(u_n)\to c>0$ and $I'(u_n)\to0$ while $\{u_n\}$ is bounded,
\[
 I(u_n)-\frac12\langle I'(u_n),u_n\rangle\to c>0,
\]
which is a contradiction. Hence \eqref{s1} is false. Thus there exist $\eta>0$ and $\{y_n\}\subset\mathbb R^N$ such that
\[
  \int_{B_r(y_n)}|u_n|^2dx\geq\eta>0.
\]
Set $\widetilde u_n(x)=u_n(x+y_n)$. Then $\{\widetilde u_n\}$ is still a bounded Palais--Smale sequence at the level $c$, and
\[
  \int_{B_r(0)}|\widetilde u_n|^2dx\geq\eta>0.
\]
Passing to a subsequence, $\widetilde u_n\rightharpoonup u_0$ in $H^2(\mathbb R^N)$ and, by the compact embedding on bounded balls, $\widetilde u_n\to u_0$ in $L^2(B_r(0))$. Hence $u_0\ne0$. Moreover, using \eqref{qslinear}, $\widetilde u_n\to u_0$ in
$L^p_{\rm loc}(\mathbb R^N)$ for every $p\in[2,2_*)$, and \eqref{f2}, we may pass to the limit in $\langle I'(\widetilde u_n),\phi\rangle$ for every $\phi\in C_c^\infty(\mathbb R^N)$; density then gives $I'(u_0)=0$. Consequently $\mathcal N\ne\emptyset$.
\end{proof}

\section{Proof of Theorem \ref{thm1}}
We shall also need the following Poho\v{z}aev identity. 
%A proof
%can be found in \cite[Lemma 3.2]{MR4291515}.
\begin{proposition}[{Poho\v{z}aev identity}]\label{pro2}
Let $u\in H^2(\mathbb R^N)$ be a weak solution of \eqref{eq1} under
$\left(f_{1}\right)$ and $\left(f_{2}\right)$. Then $u\in H^3_{\rm loc}(\mathbb R^N)$ and
\begin{align*}
P(u) := \int_{\mathbb{R}^N} \left[ \frac{N-4}{2} |\Delta u|^2 + \frac{N-2}{2} \left( |\nabla u|^2 + u^2|\nabla u|^2 \right) + \frac{N}{2} u^2 - N F(u) \right] dx = 0.
\end{align*}
\end{proposition}
\begin{proof}
Following \cite{ MR330754,MR0259693}, we know $u\in H_{\text{loc}}^3$.
The identity then follows directly from the Poho\v{z}aev variational principle applied to the functional $I$, as shown in \cite[Lemma 3.2]{MR4291515}.
\end{proof}

Define the mountain-pass level by
 $$c:=\inf_{\gamma\in \Gamma} \max _{t\in[0,1]}I(\gamma(t))$$
where
$$
\Gamma=\left\{\gamma \in C\left([0,1], H^{2}\left(\mathbb{R}^{N}\right)\right) : \gamma(0)=0\text{,\quad} I(\gamma(1))<0\right\}.
$$
By Lemma \ref{lem:mountain_pass_geometry}, the set $\Gamma\neq \emptyset$ and $c>0$.

\begin{proof}[\textbf{Proof of Theorem \ref{thm1}.}]
Following Jeanjean \cite{MR1430506}, we introduce the mapping
\begin{align*}
\varphi: E\rightarrow H^2(\mathbb{R}^N)\text{\quad by\quad} \varphi(s,v)(x)=v\left(x/e^s\right)\text{,}
\end{align*}
where $(s,v)\in E$. The composite functional $I\circ \varphi$ is given by
 \begin{align}\label{func}
 ( I\circ\varphi)(s,v)=&\frac{1}{2}\left[e^{(N-4)s}\rn|\Delta v|^2dx+e^{(N-2)s}\rn |\nabla v|^2dx+e^{Ns}\rn v^2dx\right]\nonumber\\
  &+\frac{e^{(N-2)s}}{2}\rn v^2|\nabla v|^2dx-e^{Ns}\rn F(v)dx:=\widetilde{I}(s,v).
 \end{align}
By \eqref{func}, $(f_1)$, and $(f_2)$, the functional $I\circ \varphi$ is continuously Fr\'{e}chet-differentiable on $E$. From Lemma \ref{lem:mountain_pass_geometry}, choose $v_0\in H^2(\mathbb R^N)$ and $t_0>0$ such that $I((v_0)_{t_0})<0$; equivalently, with $s^*=\log t_0$, $(I\circ\varphi)(s^*,v_0)<0$.  Consequently,
\begin{align*}
  \widetilde{\Gamma}:=\left\{\widetilde{\gamma}\in C([0,1];E): \widetilde{\gamma}(0)=(0,0)\text{,\quad} (I\circ\varphi)(\widetilde{\gamma}(1))<0\right\}\neq\emptyset.
\end{align*}
Thus, we can define the minimax level of $I\circ\varphi$,
\begin{align*}
  \widetilde{c}=\inf\limits_{\widetilde{\gamma}\in\widetilde{\Gamma}}\sup\limits_{t\in[0,1]}(I\circ \varphi)(\widetilde{\gamma}(t)).
\end{align*}
As $\Gamma=\{\varphi\circ\widetilde{\gamma}: \widetilde{\gamma}\in\widetilde{\Gamma}\}$, we have $c\leq \widetilde{c}$. For every $\gamma\in\Gamma$, the path $t\mapsto(0,\gamma(t))$ belongs to $\widetilde\Gamma$, and then $\widetilde{c}\leq c$. Hence $c=\widetilde c$. By the definition of $c$, for any $n\in\mathbb{N}$ there exists $\gamma_n\in \Gamma$, and we set $\widetilde{\gamma}_n(t)=(0,\gamma_n(t))$ such that
 $$
 \max_{t\in[0,1]}( I\circ\varphi)(\widetilde{\gamma}_n(t))=\max_{t\in[0,1]}I(\gamma_n(t))\leq c+\frac{1}{n}.
 $$
By the minimax principle \cite[Theorem 2.8]{MR1400007}, there exists a sequence $\{(s_n, v_n)\}$ in $E$ such that
\begin{enumerate}[$(a)$]
  \item $(I\circ \varphi)(s_n,v_n)\rightarrow c$,\\
  \item $\|(I\circ \varphi)'(s_n,v_n)\|_{E'}\rightarrow0 $,\\
  \item $\min_{t\in[0,1]}\|(s_n,v_n)-\widetilde{\gamma}_n(t)\|_{E}\rightarrow 0$.
\end{enumerate}
Condition $(c)$ implies that $s_n\rightarrow 0$. Setting $u_n=\varphi(s_n,v_n)$, condition $(a)$ yields $I(u_n)\rightarrow c$. For any $(\mu,\omega)\in E$, we observe
\begin{eqnarray}\label{s2}
   (I\circ \varphi)'(s_n,v_n)[\mu,\omega]=I^\prime(\varphi(s_n,v_n))[\varphi(s_n,\omega)]+P(\varphi(s_n,v_n))\mu\rightarrow0.
\end{eqnarray}
Taking $(\mu,\omega)=(1,0)$ gives $P(u_n)\to0$. To prove $I'(u_n)\to0$, let $\psi\in H^2(\mathbb R^N)$ and choose
$\omega_n(x)=\psi(e^{s_n}x)$; then $\varphi(s_n,\omega_n)=\psi$. Since $s_n\to0$, the norms of the scaling operators $\psi\mapsto\omega_n$ are uniformly bounded. Hence \eqref{s2} with $(\mu,\omega)=(0,\omega_n)$ implies $I'(u_n)\to0$ in $(H^2)'$. Therefore
\begin{align}\label{pps}
 I(u_n)\rightarrow c,\qquad I'(u_n)\rightarrow 0,\qquad P(u_n)\rightarrow0 .
\end{align}
Note that
\begin{align*}
  c+o(1)=&I(u_n)-\frac{1}{N}P(u_n)\notag\\
=&\frac{2}{N}\rn|\Delta u_n|^2dx+\frac{1}{N}\rn \left(|\nabla u_n|^2+u_n^2|\nabla u_n|^2\right)dx\notag\\
\geq& C\rn\left(|\Delta u_n|^2+|\nabla u_n|^2\right)dx.
\end{align*}
Then
\begin{align}\label{t1}
\rn|\Delta u_n|^2dx\text{,\quad }\rn |\nabla u_n|^2dx
\end{align}
 are bounded. Recall the Sobolev constant
\[S:=\inf_{\substack{u\in D^{2,2}(\mathbb{R}^N),~~
\|u\|_{2_*}=1}}\rn |\Delta u|^2dx.\]
This implies that $\|u_n\|_{2_*}$ is bounded by \eqref{t1}. By \eqref{f1}, $P(u_n)\rightarrow 0$ in \eqref{pps} and the boundedness of $\|u_n\|_{2_*}$, we obtain
\begin{align}\label{t2}
& \rn\left(\frac{N}{4} u_n^2+CN|u_n|^{2_*}\right)dx+o(1)\notag\\
\geq &N\rn F(u_n)dx+o(1)\notag\\
\geq &\frac{N}{2}\rn u_n^2dx.
\end{align}
So, $\|u_n\|_2$ is bounded. Therefore, $\{u_n\}$ is bounded in $H^2(\mathbb{R}^N)$. By Lemma \ref{lem:nonvanishing}, we obtain $\mathcal{N}\neq\emptyset$.

 Finally, we establish the existence of a ground state solution of \eqref{eq1}. Let \[m_{gs}:=\inf_{u\in\mathcal N}I(u).\]
For every $u\in\mathcal N$, Proposition \ref{pro2} gives
\[
I(u)=I(u)-\frac1N P(u)=\frac2N\int_{\mathbb R^N}|\Delta u|^2dx+\frac1N\int_{\mathbb R^N}\big(|\nabla u|^2+u^2|\nabla u|^2\big)dx>0,
\]
and hence $0\le m_{gs}<\infty$.
To complete the proof, it remains to show that $m_{gs}$ is attained on $\mathcal{N}$. Let $\{u_n\}\subset\mathcal{N}$ be such that $I(u_n)\rightarrow m_{gs}$. Since $I^\prime(u_n)=0$, Proposition \ref{pro2} gives $P(u_n)=0$. The same estimates used above show that $\{u_n\}$ is bounded in $H^2(\mathbb{R}^N)$.
Using $\langle I^\prime(u_n),u_n\rangle=0$, \eqref{f1}, and the Sobolev inequality, we obtain
\begin{align}\label{bd}
  \|u_n\|^2\leq\rn f(u_n)u_ndx\leq\frac{1}{2}\|u_n\|^2+C\|u_n\|^{2_*}.
\end{align}
From \eqref{bd}, there exists a constant $\alpha>0$ such that
\begin{align}\label{bd2}
 \liminf_{n\rightarrow\infty}\|u_n\|\geq\alpha>0.
\end{align}
If $\{u_n\}$ vanishes, i.e., if for some $r>0$,
\[
\lim_{n\to\infty}\sup_{y\in\mathbb R^N}\int_{B_r(y)}|u_n|^2dx=0,
\]
then $u_n\to0$ in $L^t(\mathbb R^N)$ for every $t\in(2,2_*)$ by Lions' lemma. Using \eqref{f2}, the boundedness of $\{u_n\}$ in $L^{2_*}(\mathbb R^N)$ and the arbitrariness of $\varepsilon$, we get
\[
 \int_{\mathbb R^N}f(u_n)u_n\,dx\to0.
\]
However, since $I'(u_n)=0$,
\[
 \|u_n\|^2+2\int_{\mathbb R^N}u_n^2|\nabla u_n|^2dx
 =\int_{\mathbb R^N}f(u_n)u_n\,dx\to0,
\]
which contradicts \eqref{bd2}. As in Lemma \ref{lem:nonvanishing}, after suitable translations there exists $u\ne0$ such that $I'(u)=0$. Thus $u\in\mathcal N$, and the definition of $m_{gs}$ gives $I(u)\ge m_{gs}$. It follows from $I'(u)=0$, $I'(u_n)=0$, weak lower semicontinuity of the norm and \eqref{qslinear} that
\begin{align*}
m_{gs}\leq I(u)&=I(u)-\frac{1}{N}P(u)\\
&=\frac{2}{N}\rn|\Delta u|^2dx+\frac{1}{N}\rn \left(|\nabla u|^2+u^2|\nabla u|^2\right)dx\\
&\leq\liminf_{n\rightarrow\infty}\left\{\frac{2}{N}\rn|\Delta u_n|^2dx+\frac{1}{N}\rn \left(|\nabla u_n|^2+ u_n^2|\nabla u_n|^2\right)dx\right\}\\
&\leq\liminf_{n\rightarrow\infty}\left(I(u_n)-\frac{1}{N}P(u_n)\right)=m_{gs}.
\end{align*}
Thus, $I(u)=m_{gs}$, which proves that $u$ is a ground state solution.
\end{proof}

\section{Proof of Theorem \ref{th2}}

 To construct the required minimax levels, we use the augmented-dimension trick introduced in \cite{Hirata2010}. Choose once and for all an exponent
\[
       p_0\in(3,2_*-1).
\]
For $s>0$, define
\[
 m(s)=\sup_{0<\tau\le s}\frac{\hat{f}(\tau)}{\tau^{p_0}},
\]
where
\begin{align*}
       \hat{f}(s)= \begin{cases}
        \max\left\{f(s)-\frac{s}{2},0\right\}, & \text{for } s\geq 0,\\
        - \hat{f}(-s), & \text{for } s<0.
       \end{cases}
\end{align*}
 Define
\begin{align*}
      g(s)=
     \begin{cases} s^{p_0}m(s), & \text{for } s>0,\\ 0, & \text{for } s=0,\\ -g(-s), & \text{for } s<0,
     \end{cases}
\end{align*}
and set $G(s)=\int_0^s g(t)\,dt$.

We introduce the comparison functional $J: H_r^2(\mathbb{R}^N)\rightarrow\mathbb{R}$ by
\begin{align}\label{cf}
  J(u)=\frac{1}{2} \int_{\mathbb{R}^{N}}(|\Delta u|^2+|\nabla u|^{2}+u^2|\nabla u|^2)\, dx
       +\frac14\int_{\mathbb R^N}u^2\,dx
       - \int_{\mathbb{R}^{N}}G(u)\, dx,
\end{align}
for $u\in H_r^2(\mathbb{R}^N)$. Since $g(s)\ge (f(s)-s/2)^+$ for $s\ge0$ and $f$ is odd, we have
\[
   G(s)\ge F(s)-\frac{s^2}{4}\qquad(s\in\mathbb R).
\]
Consequently,
\[
 I(u)-J(u)=\int_{\mathbb R^N}\left(G(u)-F(u)+\frac14u^2\right)dx\ge0,
\]
that is, $I(u)\ge J(u)$ for all $u\in H_r^2(\mathbb R^N)$. Following the reasoning in \cite[Lemma 2.1, Corollary 2.2, Lemma 2.4]{Hirata2010}, we obtain the following lemma:
 \begin{lemma}[\cite{Hirata2010}]\label{lem41}
 \begin{enumerate}[$(i)$]
   \item  $0\leq (p_0+1)G(s)\leq s g(s)$,
   \item $\lim_{|s|\rightarrow+\infty}\frac{g(s)}{|s|^{2_*-1}}=0$,

   \item There exists $\delta_0>0$ such that $g(s)=0$ for $|s|\leq\delta_0$,
   \item $J\in C^1(H_r^2(\mathbb{R}^N))$, and, for all $u,\varphi\in H_r^2(\mathbb{R}^N)$,
   \begin{align*}
   \langle J'(u), \varphi\rangle =& \int_{\mathbb{R}^{N}} \left( \Delta u \Delta \varphi + \nabla u \cdot \nabla \varphi + \frac12 u \varphi + u\varphi|\nabla u|^2 + u^2 \nabla u \cdot \nabla \varphi \right) dx\notag\\
&\quad- \int_{\mathbb{R}^{N}} g(u)\varphi dx .
\end{align*}
   %\item For $N\geq 2$, there exists a $C_N>0$ such that $|u(x)|\leq C_N|x|^{-\frac{N-1}{2}}\|u\|_{H^1}$ for $u\in H_r^1$ and $|x|\geq 1$.
 \end{enumerate}

 \end{lemma}

%\begin{proof}
 % Since the mapping
%  \begin{align*}
 %   (0, \infty) \to \mathbb{R};\quad s \mapsto g(s) / s^{p_0}=m(s)
%  \end{align*}
%  is non-decreasing, we have for \( s> 0 \)

%\begin{align*}
%sg(s) - (p_0 + 1) G(s)
%=& \int_0^s\left(g(s) - (p_0 + 1) g(\tau)\right) \, d\tau\\
%=& \int_0^s \left(s^{p_0} \frac{g(s)}{s^{p_0}} - (p_0 + 1) \tau^{p_0} \frac{g(\tau)}{\tau^{p_0}}\right) \, d\tau\\
%\geq& \int_0^s \left(s^{p_0} \frac{g(s)}{s^{p_0}} - (p_0 + 1) \tau^{p_0} \frac{g(s)}{s^{p_0}}\right) \, d\tau = 0.
%\end{align*}
%Therefore (i) holds.
%\end{proof}

\begin{definition}
%\label{def:pi_k_minus_1}
For any integer $k \geq 1$, we define the $(k-1)$-dimensional polyhedron:
\begin{equation*}
\pi_{k-1} =\left\{\ell=(\ell_1,\ldots,\ell_k) \in \mathbb{R}^k \ \middle| \ \sum_{j=1}^k |\ell_j| = 1 \right\}
\end{equation*}
which is the $\ell^1$-unit sphere in $\mathbb{R}^k$.
\end{definition}

\begin{remark}
The polyhedron $\pi_{k-1}$ possesses exactly $2k$ vertices located at $\pm \mathbf{e}_i$ (the standard basis vectors) and $2^k$ facets.
Since \( \pi_{k-1} \) is homeomorphic to \( S^{k-1} \) by an odd homeomorphism (take for instance the radial projection \( \pi_{k-1} \to S^{k-1} \)),
 it is clear that \(\operatorname{genus}(\pi_{k-1}) = k \), where \(\operatorname{genus}(Y)\) is the Krasnosel'skii genus of \(Y\), and
\[
\operatorname{genus}(Y) := \min \left\{ n \in \mathbb{N}^+ \mid \exists \, \varphi : Y \to \mathbb{R}^n \setminus \{0\}, \ \varphi \text{ is odd and continuous} \right\}.
\]
We set $\operatorname{genus}(Y) = \infty$ if such a map $\varphi$ does not exist, and set $\operatorname{genus}(Y) = 0$ if $Y = \emptyset$.
 %We will prove that \( \beta_k > 0 \) using the next result.
%
%Through radial projection, it is homeomorphic to $S^{k-1}$, which immediately implies that $\gamma(\pi_{k-1})=k$, where $\gamma$ denotes the Krasnosel'skii genus.
\end{remark}

\begin{lemma}
%\label{lem5}
There exists a function $\phi\in C^2(\mathbb{R})$ such that
\[
\phi(t)=
\begin{cases}
0, & t\le 0,\\
10t^3-15t^4+6t^5, & 0<t<1,\\
1, & t\ge 1.
\end{cases}
\]
Moreover,
\[
0\le \phi\le 1,\qquad
\phi(0)=0,\quad \phi(1)=1,
\]
\[
\phi'(0)=\phi'(1)=0,\qquad
\phi''(0)=\phi''(1)=0,
\]
and there exists an absolute constant $C_\phi>0$ such that
\[
    \|\phi'\|_{L^\infty(\mathbb{R})}+\|\phi''\|_{L^\infty(\mathbb{R})}\le C_\phi .
\]
\end{lemma}

\begin{proof}
For $0<t<1$, we set
\[
    q(t)=10t^3-15t^4+6t^5.
\]
Then
\[
    q'(t)=30t^2(1-t)^2,
    \qquad
    q''(t)=60t(1-t)(1-2t).
\]
Thus $q(0)=0$, $q(1)=1$, and
\[
q'(0)=q'(1)=q''(0)=q''(1)=0.
\]
Consequently, the extension defined above by constants on $(-\infty,0]$ and $[1,\infty)$
belongs to $C^2(\mathbb{R})$. The stated bounds follow from the explicit formulas for $q'$ and $q''$.
\end{proof}

\begin{proposition}
\label{prop:BL-nodal-map}
Let $k\ge 1$. Assume that $f$ is odd and $(f_3)$ holds, then there exists $R=R(k)>0$ and an odd continuous mapping
\[
    \tau:\pi_{k-1}\longrightarrow H_{0,r}^2(B_R),
\]
such that, for every $\ell\in\pi_{k-1}$, the function $u=\tau(\ell)$ satisfies
\begin{enumerate}[$(i)$]
    \item $0\notin \tau(\pi_{k-1})$;
    \item there exist constants $\varrho_k,C_k>0$, depending on $k$ but independent of $\ell$, such that
    \[
        \varrho_k\le \|\Delta u\|_{L^2(B_R)}\le C_k;
    \]
    \item
    \[
        \int_{B_R}\widetilde{F}(u)d x\ge 1.
    \]
\end{enumerate}
\end{proposition}

\begin{proof}
Since $f$ is odd, $F$ and $\widetilde{F}$ are even. Hence
\[
    a:=\widetilde{F}(\xi)=\widetilde{F}(-\xi)>0.
\]
Set
\[
    M:=\max_{|s|\le \xi}|\widetilde{F}(s)|.
\]

\medskip
\noindent\textbf{Step 1. A fixed-width smoothing procedure.}
Fix $h=1$. Let $\varrho\in C_c^\infty((-1,1))$ be even, nonnegative, and satisfy
\[
    \int_{\mathbb{R}} \varrho(t) d t=1,
\]
and set $\varrho_h(t)=h^{-1}\varrho(t/h)$. We shall choose $R>4h$ sufficiently large later.
Define the boundary cut-off function
\[
    \chi_R(r):=\phi\left(\frac{r-h}{h}\right)
             \phi\left(\frac{R-h-r}{h}\right),
    \qquad 0\le r\le R.
\]
Then $\chi_R\in C^2([0,R])$, $0\le \chi_R\le 1$,
\[
    \chi_R(r)=0\quad\text{for }0\le r\le h\text{ and }R-h\le r\le R,
\]
while
\[
    \chi_R(r)=1\quad\text{for }2h\le r\le R-2h.
\]
Moreover, $\chi_R'(0)=\chi_R'(R)=0$, and
\[
    \|\chi_R\|_{C^2([0,R])}\le C_h,
\]
where $C_h$ is independent of $R$.

For $\ell=(\ell_1,\ldots,\ell_k)\in\pi_{k-1}$, define
\[
    \alpha_i(\ell):=R\sum_{j=1}^i |\ell_j|,
    \qquad i=0,1,\ldots,k,
\]
with $\alpha_0=0$ and $\alpha_k=R$. Define a step function $s_\ell$ on $[0,R]$ by
\[
    s_\ell(r) = \operatorname{sgn}(\ell_i)
    \quad\text{for } r\in (\alpha_{i-1}(\ell),\alpha_i(\ell))\text{ and }\ell_i\ne 0.
\]
Intervals corresponding to $\ell_i=0$ have length zero and are irrelevant. We extend
$s_\ell$ by zero outside $[0,R]$ and set
\[
    v_\ell:=\varrho_h*s_\ell\quad\text{on }\mathbb{R}.
\]
Finally, define the radial function
\[
    \tau(\ell)(x):=\xi\,\chi_R(|x|)v_\ell(|x|),
    \qquad x\in B_R.
\]
Since $\chi_R$ vanishes in neighborhoods of both endpoints $0$ and $R$, the function
$x\mapsto\tau(\ell)(x)$ is $C^2$ in $\overline{B_R}$, vanishes near $\partial B_R$,
and is radial. Therefore
\[
    \tau(\ell)\in H_{0,r}^2(B_R).
\]
The construction is odd because $s_{-\ell}=-s_\ell$, hence
\[
    \tau(-\ell)=-\tau(\ell).
\]

\medskip
\noindent\textbf{Step 2. Continuity of $\tau$.}
We claim that $\ell\mapsto s_\ell$ is continuous from $\pi_{k-1}$ into $L^1(\mathbb{R})$. Indeed, if $\ell^m\to\ell$ in $\pi_{k-1}$, then the endpoints $\alpha_i(\ell^m)$ converge to $\alpha_i(\ell)$ for every $i$. For every index with $\ell_i\ne0$, the signs of $\ell_i^m$ and $\ell_i$ agree for all sufficiently large $m$; the corresponding intervals differ only through their endpoints. Indices with $\ell_i=0$ contribute intervals of length $R|\ell_i^m|$. Consequently,
\[
 \|s_{\ell^m}-s_\ell\|_{L^1(\mathbb R)}
 \le 2\sum_{i=1}^{k-1}|\alpha_i(\ell^m)-\alpha_i(\ell)|
 +2R\sum_{\{i:\ell_i=0\}}|\ell_i^m|\longrightarrow0.
\]
Since convolution with the fixed smooth kernel $\varrho_h$ maps $L^1(\mathbb{R})$ continuously into $C^2_{\mathrm{loc}}(\mathbb{R})$, we have
\[
    \|v_{\ell^m}-v_\ell\|_{C^2([0,R])}
    \le C_h\|s_{\ell^m}-s_\ell\|_{L^1(\mathbb{R})}\to 0.
\]
Multiplication by the fixed cut-off $\chi_R$ then yields
\[
    \tau(\ell^m)\to \tau(\ell)\quad\text{in }H^2(B_R).
\]
Thus $\tau:\pi_{k-1}\to H_{0,r}^2(B_R)$ is continuous.

\medskip
\noindent\textbf{Step 3. The nonlinear lower bound.}
Let
\[
    \mathcal I_\ell:=\{\alpha_i(\ell):i=1,\cdots,k-1\}
\]
be the set of interior interfaces. Define the exceptional radial set
\[
   D_\ell := \left\{ r \in [0, R] : \operatorname{dist}(r, \{0, R\} \cup \mathcal{I}_\ell) \leq 2h \right\}.
\]
Outside $D_\ell$, the support of the mollifier does not intersect any interface or endpoint and
$\chi_R=1$. Therefore
\[
    \tau(\ell)(x)=\pm\xi
    \quad\text{whenever } |x|\notin D_\ell.
\]
The set $D_\ell$ is the union of at most $k+1$ intervals of length at most $4h$.
Consequently, for $R\ge 1$,
\[
    \left|\{x\in B_R: |x|\in D_\ell\}\right|
    \le C_N(k+1)hR^{N-1},
\]
where $C_N>0$ depends only on $N$. Since $|\tau(\ell)|\le \xi$, we obtain
\begin{align*}
    \int_{B_R}\widetilde{F}(\tau(\ell))d x
    &\ge a\,|B_R|- (a+M)C_N(k+1)hR^{N-1}.
\end{align*}
Choosing $R=R(k)>4h$ sufficiently large gives
\[
    \int_{B_R}\widetilde{F}(\tau(\ell)) d x\ge 1
    \quad\text{for all }\ell\in\pi_{k-1}.
\]
This proves item \textup{(iii)} and also implies $0\notin\tau(\pi_{k-1})$.

\medskip
\noindent\textbf{Step 4. Uniform upper and lower bounds for $\|\Delta\tau(\ell)\|_2$.}
The functions $s_\ell$ satisfy $|s_\ell|\le 1$. Hence $v_\ell$ and its first two
derivatives are uniformly bounded in terms of $h$ and $\varrho$ only. Together with the
uniform $C^2$ bound for $\chi_R$, this yields
\[
    \|\tau(\ell)\|_{C^2([0,R])}\le C(\xi,h,\varrho)
    \quad\text{for all }\ell\in\pi_{k-1}.
\]
Moreover $\tau(\ell)$ vanishes on $[0,h]$. Therefore, for $r\ge h$,
\[
    |\Delta \tau(\ell)(r)|
    \le |\tau(\ell)''(r)|+\frac{N-1}{r}|\tau(\ell)'(r)|
    \le C(\xi,h,\varrho,N).
\]
It follows that
\[
   \|\Delta\tau(\ell)\|_{L^2(B_R)}\le C_k
\]
for a constant $C_k>0$ independent of $\ell$; here $C_k$ may depend on $R(k)$.

Finally, since $\pi_{k-1}$ is compact and $\tau$ is continuous, the set
$\tau(\pi_{k-1})$ is compact in $H_{0,r}^2(B_R)$. By item \textup{(iii)}, it does not
contain the zero function. Since $\|\Delta\cdot\|_{L^2(B_R)}$ is an equivalent norm
on $H_{0,r}^2(B_R)$, the continuous function
\[
    \ell\mapsto \|\Delta\tau(\ell)\|_{L^2(B_R)}
\]
has a strictly positive minimum on $\pi_{k-1}$. Thus there exists $\varrho_k>0$ such
that
\[
    \varrho_k\le \|\Delta\tau(\ell)\|_{L^2(B_R)}\le C_k
    \quad\text{for all }\ell\in\pi_{k-1}.
\]
The proof is complete.
\end{proof}

The functionals $I$ and $J$ satisfy the following geometric properties.

\begin{proposition}\label{p1}
We have the following properties:
\begin{enumerate}[$(1)$]
    \item  There exist $\rho,\mu>0$ such that
     \begin{align*}
       I(u) \geq J(u) \geq 0\text{\quad for all~~}\|u\| \leq \rho,\\
       I(u) \geq J(u) \geq \mu\text{\quad for all~~}\|u\| = \rho;
     \end{align*}

    \item  For any $n\in\mathbb{N}^+$, there is an odd continuous mapping
\[
\gamma_{0n} : S^{n-1} \to H_r^2(\mathbb{R}^N)
\]
where $S^{n-1} = \{\sigma \in \mathbb{R}^n \mid |\sigma| = 1\}$, such that

\begin{align}
\label{y2}
J(\gamma_{0n}(\sigma)) \leq I(\gamma_{0n}(\sigma)) < 0\text{,\quad for all }\sigma \in S^{n-1}.
\end{align}
\end{enumerate}
\end{proposition}

\begin{proof}
As in the proof of Lemma \ref{lem:mountain_pass_geometry}, by $(ii)$ and $(iii)$ of Lemma \ref{lem41} there exist $\rho,\mu>0$ such that
\[
     J(u)\ge\mu>0\quad\hbox{if }\|u\|=\rho,
     \qquad J(u)\ge0\quad\hbox{if }\|u\|\le\rho.
\]
Together with $I\ge J$, this proves item $(1)$.

We now prove item $(2)$. By Proposition \ref{prop:BL-nodal-map}, applied with $k=n$, and by a fixed odd homeomorphism $h_n:S^{n-1}\to\pi_{n-1}$, there exists an odd continuous map
\[
       \zeta_n:S^{n-1}\to H_{0,r}^2(B_R)\subset H_r^2(\mathbb R^N)
\]
such that $0\notin \zeta_n(S^{n-1})$ and
\[
       \int_{\mathbb R^N}\widetilde F(\zeta_n(\sigma))dx\ge1
       \quad\hbox{for all }\sigma\in S^{n-1}.
\]
For $\theta\ge1$, define
\[
      \gamma_{0n}(\sigma)(x)=\zeta_n(\sigma)(x/\theta),
      \qquad \sigma\in S^{n-1}.
\]
Then $\gamma_{0n}$ is odd and continuous. Since
\begin{align*}
  I(\gamma_{0n}(\sigma))
  &=\frac12\theta^{N-4}\int_{\mathbb R^N}|\Delta\zeta_n(\sigma)|^2dx
    +\frac12\theta^{N-2}\int_{\mathbb R^N}\left(|\nabla\zeta_n(\sigma)|^2+\zeta_n(\sigma)^2|\nabla\zeta_n(\sigma)|^2\right)dx\\
  &\quad -\theta^N\int_{\mathbb R^N}\widetilde F(\zeta_n(\sigma))dx,
\end{align*}
and the positive terms above are uniformly bounded in $\sigma$, choosing $\theta=\theta_n$ sufficiently large gives
\[
        I(\gamma_{0n}(\sigma))<0\quad\hbox{for all }\sigma\in S^{n-1}.
\]
Because $J\le I$, we also have $J(\gamma_{0n}(\sigma))\le I(\gamma_{0n}(\sigma))<0$. This proves item $(2)$.
\end{proof}

\begin{lemma}\label{J1}
  The functional $J|_{H_r^2}$ satisfies the Palais--Smale condition.
\end{lemma}

\begin{proof}
Let $\{u_m\}\subset H_r^2(\mathbb R^N)$ be a Palais--Smale sequence for $J$ at level $c$, that is
\[J(u_m)\rightarrow c\text{\quad and \quad} \|J|_{H_r^2}'(u_m)\|_{(H_r^2)'}\rightarrow 0\text{.}\]
 Put
\[
 A_m=\int_{\mathbb R^N}(|\Delta u_m|^2+|\nabla u_m|^2)dx,
 \quad U_m=\int_{\mathbb R^N}u_m^2dx.
\]
Using Lemma \ref{lem41} $(i)$ and $p_0\in(3, 2_*-1)$, we obtain
\begin{align*}
 c+o(1)(1+\|u_m\|)
 &\geq J(u_m)-\frac{1}{p_0+1}\langle J'(u_m),u_m\rangle\\
 &=\left(\frac12-\frac{1}{p_0+1}\right)A_m
   +\left(\frac14-\frac{1}{2(p_0+1)}\right)U_m
   +\left(\frac12-\frac{2}{p_0+1}\right)Q(u_m)\\
 &\quad +\int_{\mathbb R^N}\left(\frac{1}{p_0+1}g(u_m)u_m-G(u_m)\right)dx\\
 &\ge C\bigl(A_m+U_m+Q(u_m)\bigr).
\end{align*}
Hence $\{u_m\}$ is bounded in $H_r^2(\mathbb R^N)$. Up to a subsequence,
\[
 u_m\rightharpoonup u\quad\hbox{in }H_r^2(\mathbb R^N),\qquad
 u_m\to u\quad\hbox{in }L^q(\mathbb R^N),\ 2<q<2_*,
\]
and $u_m(x)\to u(x)$ a.e. in $\mathbb R^N$. Since $g$ satisfies $(ii)$ and $(iii)$ of Lemma \ref{lem41}, for every $\varepsilon>0$ there exists $C_{\varepsilon}>0$ such that
\begin{align}\label{g12}
|g(s)|\leq \varepsilon\bigl(|s|^{2_*-1}+|s|\bigr)+C_{\varepsilon}|s|^p\text{\quad for~~}p\in(1,2_*-1).
\end{align}
By \eqref{g12}, the H\"{o}lder inequality and $u_m\rightarrow u$ in $L^q$ with $q\in (2,2_*)$ we obtain
\begin{align*}
     \left| \int_{\mathbb R^N}g(u_m)(u_m-u)dx\right|\leq&\varepsilon\|u_m\|_2\|u_m-u\|_2+\varepsilon\|u_m\|_{2_*}^{2_*-1}\|u_m-u\|_{2_*}\notag\\
      &\quad+C_\varepsilon\|u_m\|_{p+1}^{p}\|u_m-u\|_{p+1}\notag\\
      \leq& C\varepsilon+o(1) .
\end{align*}
Since $\varepsilon>0$ is arbitrary,
\begin{equation*}
 \int_{\mathbb R^N}g(u_m)(u_m-u)\,dx=o(1).
\end{equation*}
Since $u_m\rightharpoonup u$ in $H_r^2(\mathbb{R}^N)$, one can show in a standard way that
\[
 \int_{\mathbb R^N}(g(u_m)-g(u))u\,dx\to0.
\]
Consequently,
\begin{align}\label{g111}
 \int_{\mathbb R^N}g(u_m)u_m\,dx\to\int_{\mathbb R^N}g(u)u\,dx.
\end{align}
Passing to the limit in $\langle J|_{H_r^2}'(u_m),\phi\rangle$ for every radial $\phi\in C_c^\infty(\mathbb R^N)$ and then using density gives $J|_{H_r^2}'(u)=0$. Hence, by \eqref{g111} and $J|_{H_r^2}'(u)=0$ we have
\begin{align*}
 A+\frac{1}{2}U+2Q(u)
 &=\int_{\mathbb R^N}g(u)u\,dx\\
 &=\lim_{m\to\infty}\left(A_m+\frac12U_m+2Q(u_m)\right),
\end{align*}
where $A$ and $U$ are the corresponding quantities for $u$. Since the quadratic part and $Q$ are weakly lower semicontinuous by \eqref{qslinear}, equality of the sums forces
\[
 A_m+\frac12U_m\to A+\frac12U,
 \qquad Q(u_m)\to Q(u).
\]
The weak convergence and convergence of the quadratic norm imply $u_m\to u$ strongly in $H_r^2(\mathbb R^N)$. Thus $J|_{H_r^2}$ satisfies the Palais--Smale condition.
\end{proof}

By Proposition \ref{p1}, both functionals $I$ and $J$ exhibit the symmetric mountain-pass geometry, allowing us to define corresponding minimax values. Following the framework of \cite[Chapter 9]{Rabinowitz1986}, we define for each $n\in\mathbb{N}^+$ the minimax values
\begin{equation}\label{eq:3.1}
b_{n}=\inf_{\gamma\in\Gamma_{n}}\max_{\sigma\in D_{n}}I(\gamma(\sigma))\text{,\quad}
c_{n}=\inf_{\gamma\in\Gamma_{n}}\max_{\sigma\in D_{n}}J(\gamma(\sigma))\text{,}
\end{equation}
where \begin{itemize}
    \item $D_n = \{\sigma = (\sigma_1,\ldots,\sigma_n) \in \mathbb{R}^n : |\sigma| \leq 1\}$.
    \item $\Gamma_n$ denotes the class of admissible mappings defined by
    \begin{equation}\label{eq:3.2}
    \Gamma_n = \left\{\gamma \in C(D_n, H^2_r(\mathbb{R}^N)) \colon
    \begin{aligned}
        &\gamma(-\sigma) = -\gamma(\sigma) \text{ for all } \sigma \in D_n, \\
        &\gamma(\sigma) = \gamma_{0n}(\sigma) \text{ for all } \sigma \in \partial D_n
    \end{aligned}
    \right\},
    \end{equation}
    where $\gamma_{0n}: \partial D_n = S^{n-1} \to H^2_r(\mathbb{R}^N)$ is the mapping constructed in Proposition \ref{p1}.
\end{itemize}
We note that $\Gamma_n$ is nonempty for all $n$, as it contains the explicit mapping:
\[
\gamma(\sigma)=\left\{\begin{array}{ll}
|\sigma|\gamma_{0n}\left(\frac{\sigma}{|\sigma|}\right) & \text{for }\sigma\in D_{n}\setminus\{0\},\\
0 & \text{for }\sigma=0\text{,}
\end{array}\right.
\]
which satisfies all required conditions.

Next, for each $j\in\mathbb{N}^+$, define the minimax families $\Lambda_j$ by
\[
\Lambda_j = \{h(\overline{D_m \setminus Y}) : h \in \Gamma_m, \, m \geq j, \, Y \in \mathcal{E}_m \text{ and } \operatorname{genus}(Y) \leq m - j\},
\]
where
\[
\mathcal E_m:=\{Y\subset \mathbb{R}^m\setminus\{0\}:Y\text{ is closed and symmetric}\}.
\]
We also write $\mathcal E(H_r^2)$ for the family of closed symmetric subsets of
$H_r^2(\mathbb R^N)\setminus\{0\}$. This construction is adapted from
\cite[Chapter 9, p.~56]{Rabinowitz1986}.

The next proposition records the basic properties of $\Lambda_j$.
\begin{proposition}\label{p3}
The sets $\Lambda_j$ possess the following properties:
\begin{enumerate}[$(i)$]
 \item $\Lambda_j\neq \emptyset$ for all $j\in \mathbb{N}^+$;
 \item (Monotonicity) $\Lambda_{j+1}\subset\Lambda_j$;
 \item (Invariance) If $\varphi\in C(H_r^2,H_r^2)$ is odd and
 $\varphi=\mathrm{id}$ on $\gamma_{0m}(\partial D_m)$ for every $m\ge j$,
 then $\varphi(B)\in\Lambda_j$ for every $B\in\Lambda_j$;
 \item (Excision) If $B\in\Lambda_j$ and
 $Z\in\mathcal E(H_r^2)$ satisfies $\operatorname{genus}(Z)\le i<j$, then $\overline{B\backslash Z}\in \Lambda_{j-i}$.
\end{enumerate}
\end{proposition}

\begin{proof}
$(i)$ For every $m$, the cone extension
\[
 h_m(\sigma)=
 \begin{cases}
 |\sigma|\gamma_{0m}(\sigma/|\sigma|),&\sigma\ne0,\\
 0,&\sigma=0,
 \end{cases}
\]
belongs to $\Gamma_m$. Taking $m=j$ and $Y=\emptyset$ shows that
$\Lambda_j\ne\emptyset$.

$(ii)$ If $B = h(\overline{D_m \setminus Y}) \in \Lambda_{j+1}$, by the definition of $\Lambda_{j}$ we have
  \[m \geq j+1 \geq j\text{,\quad} h \in \Gamma_m\text{,\quad } Y \in \mathcal{E}_m\]
 and $\operatorname{genus}(Y)  \leq m - (j+1) \leq m - j$. Therefore $B \in \Lambda_j$, that is, $\Lambda_{j+1}\subset\Lambda_j$.

$(iii)$ If $B=h(\overline{D_m\setminus Y})\in\Lambda_j$, then
$\varphi\circ h$ is odd and continuous. On $\partial D_m$,
\[
 (\varphi\circ h)(\sigma)=\varphi(\gamma_{0m}(\sigma))
 =\gamma_{0m}(\sigma),
\]
so $\varphi\circ h\in\Gamma_m$ and
$\varphi(B)=(\varphi\circ h)(\overline{D_m\setminus Y})\in\Lambda_j$.

$(iv)$ Let $B = h(\overline{D_m \setminus Y}) \in \Lambda_j$ and $Z\in\mathcal E(H_r^2)$ with $\operatorname{genus}(Z) \leq i < j$. We claim
\begin{equation}\label{excision}
\overline{B \setminus Z} = h\bigl(\overline{D_m \setminus (Y \cup h^{-1}(Z))}\bigr).
\end{equation}
 Since $h(0)=0$ and $0\notin Z$, the set $h^{-1}(Z)$ is a closed
symmetric subset of $D_m\setminus\{0\}$. Moreover, by monotonicity of the genus
under odd continuous maps,
\[
 \operatorname{genus}(h^{-1}(Z))\le\operatorname{genus}(Z)\leq i<j.
\]
Thus
\begin{align*}
 \operatorname{genus}\bigl(Y\cup h^{-1}(Z)\bigr)
 &\le\operatorname{genus}(Y)+\operatorname{genus}(h^{-1}(Z))\\
 &\le m-j+i=m-(j-i),
\end{align*}
which proves $\overline{B \setminus Z }\in \Lambda_{j-i}$ thanks to \eqref{excision}.

Now, we shall prove \eqref{excision}. Let $y \in h(\overline{D_m \setminus (Y \cup h^{-1}(Z))})$. Then there exists
$x \in \overline{D_m \setminus (Y \cup h^{-1}(Z))} \subset D_m$ such that $y = h(x)$. By definition of closure,
there is a sequence $\{x_n\} \subset D_m \setminus (Y \cup h^{-1}(Z))$ such that $x_n \to x$.
By continuity of $h$, $h(x_n) \to h(x) = y$. By $x_n \in D_m \setminus (Y \cup h^{-1}(Z))$, we have \[h(x_n) \in h(D_m \setminus Y) \setminus Z \subset \overline{B \setminus Z}.\]
This implies that $y\in \overline{B \setminus Z}$. Therefore
\begin{align}
\label{exci1} h\bigl(\overline{D_m \setminus (Y \cup h^{-1}(Z))}\bigr)\subset \overline{B \setminus Z} .
\end{align}
On the other hand, for any $y \in B \setminus Z$ there exists $y = h(w)$ where
\[
w \in \overline{D_m \setminus Y} \setminus h^{-1}(Z) \subset\overline{ D_m \setminus (Y \cup h^{-1}(Z))}.
\]
Together with the continuity of $h$, we obtain
\begin{equation}\label{exci2}
\overline{B \setminus Z} \subset h(\overline{D_m \setminus (Y \cup h^{-1}(Z))}).
\end{equation}
The claim \eqref{excision} follows from \eqref{exci1} and \eqref{exci2}.

\end{proof}
We define another sequence of minimax values by
\begin{equation}\label{mini1}
  d_j = \inf_{A \in \Lambda_j} \max_{u \in A} J(u).
\end{equation}
For the radius $\rho$ fixed in Proposition \ref{p1}, set
\[
\mathbb B_\rho:=\{u\in H_r^2(\mathbb R^N):\|u\|<\rho\}.
\]
\begin{proposition}
  \label{p4}
  If $B\in \Lambda_j$, then $B\cap \left(\partial\mathbb B_\rho\cap H_r^2\right)\neq\emptyset$.
\end{proposition}
\begin{proof}
  Set $B = h(\overline{D_m \setminus Y})$ where $h\in \Gamma_m$, $m \geq j$, $Y\in\mathcal{E}_m$ and $\operatorname{genus}(Y) \leq m - j$. For any $\sigma\in \partial D_m$, by \eqref{eq:3.2} and \eqref{y2} we obtain that
  \[
    J(h(\sigma))=J(\gamma_{0m}(\sigma))<0.
  \]
Together with Proposition \ref{p1} $(1)$, we have $\|h(\sigma)\|>\rho$.
  Let
$\hat{O} = \{\sigma \in D_m \mid h(\sigma) \in \mathbb B_\rho\}$. Since $h$ is odd, $0 \in \hat{O}$.
Let $O$ denote the component of $\hat{O}$ containing $0$. Since $\hat O$ is symmetric and $0\in O$, the uniqueness of the connected component containing $0$ gives $-O=O$. Moreover, the boundary condition $\|h(\sigma)\|>\rho$ on $\partial D_m$ implies $\overline O\subset\operatorname{int}D_m$. Thus $O$ is a bounded symmetric neighborhood of $0$ in $\mathbb{R}^m$.
Therefore, $\operatorname{genus}(\partial O) = m$ and
\begin{align*}
  \operatorname{genus}(\overline{\partial O \setminus Y}) \geq & \operatorname{genus}(\partial O)-  \operatorname{genus}(Y)\\
  \geq& m - (m - j) \\
  =& j \geq 1.
\end{align*}
In particular, $\overline{\partial O \setminus Y}\neq \emptyset$. By the definition of $O$ and continuity of $h$, $h(\partial O) \subset \partial\mathbb B_\rho$. Therefore,
$h(\overline{\partial O \setminus Y})\subset (B\cap \partial\mathbb B_\rho)$, that is, $B\cap \left(\partial\mathbb B_\rho\cap H_r^2\right)\neq\emptyset$.
\end{proof}
Propositions \ref{p3} and \ref{p4} imply that $0<\mu\le d_j<+\infty$, where $\mu$ is given in Proposition \ref{p1}.

  For \( d \in \mathbb{R} \) we set
\[
K^d = \left\{ u \in H_r^2(\mathbb{R}^N) : J(u) =d \text{\quad and \quad} J|_{H_r^2}'(u) = 0 \right\}.
\]
Then \(K^d\) is compact by Lemma \ref{J1}. For any \( d \in \mathbb{R} \), we define the sublevel sets of \( J\) as follows:
\[
J^d := \{ u \in H_r^2(\mathbb{R}^N) : J(u) \leq d \}.
\]
\begin{lemma}\label{lem49}
  Let $K^d$ be the critical set defined above. If
\[
d_{j} = d_{j+1} = \cdots = d_{j+i} := d \geq \mu>0
\]
for some $j \geq 1$ and $i \geq 0$, then
\[
\operatorname{genus}(K^d) \geq i + 1.
\]
In particular, if $i \geq 1$, then the critical set at level d is infinite.
\end{lemma}

\begin{proof}
Assume by contradiction that $\operatorname{genus}(K^d)\le i$. By Lemma \ref{J1}, $K^d$ is compact. By item~$5^\circ$ of \cite[Proposition 7.5]{Rabinowitz1986}, there exists a symmetric open neighborhood $\mathcal N_d$ of $K^d$, with $0\notin\overline{\mathcal N_d}$, such that
\[
 \operatorname{genus}(\overline{\mathcal N_d})=\operatorname{genus}(K^d)\le i.
\]
Apply the deformation theorem \cite[Theorem A.4]{Rabinowitz1986} with  $\mathcal O=\overline{\mathcal N_d}$ and $\overline\varepsilon=\mu/2$. There exist $\varepsilon\in(0,\overline\varepsilon)$ and $\eta\in C([0,1]\times H_r^2,H_r^2)$ such that $\eta(1,\cdot)$ is odd,
\begin{equation}\label{y3}
 \eta(1,J^{d+\varepsilon}\setminus\mathcal{O})\subset J^{d-\varepsilon},
\end{equation}
and $\eta(1,u)=u$ for $u\in J^{d-\overline\varepsilon}$.

By the definition of $d_{j+i}=d$, choose $A\in\Lambda_{j+i}$ such that
\[\max_AJ\le d+\varepsilon,\]
that is \( A \subset J^{d+\epsilon} \). Proposition \ref{p3} $(iv)$, applied with  $Z=\mathcal{O}$, yields a set $\overline{A \setminus \mathcal{O}}\in \Lambda_j$ satisfying
\[
\overline{A \setminus \mathcal{O}}
 \subset \overline{J^{d+\varepsilon}\setminus\mathcal{O}}.
\]
Therefore \eqref{y3} gives
\begin{equation}\label{y4}
 \eta(1,\overline{A \setminus \mathcal{O}})\subset \eta(1,\overline{J^{d+\varepsilon}\setminus\mathcal{O}})\subset\overline{\eta(1,J^{d+\varepsilon}\setminus\mathcal{O})}\subset J^{d-\varepsilon}.
\end{equation}
and
\[
\max_{\eta(1, \overline{A \setminus \mathcal{O}})} J \leq d - \epsilon < d.
\]
For any $u\in \gamma_{0m}(\partial D_m)$ satisfying $J(u)<0$, we have
\[J(u)\leq 0< d-\overline{\varepsilon}.\]
Hence $2^\circ$ of \cite[Theorem A.4]{Rabinowitz1986} and our choice of $\overline{\varepsilon}$ imply $\eta(1,u)=u$ on $J(u)<d-\overline{\varepsilon}$. Consequently, by $(iii)$ of Proposition \ref{p3} we have $\eta(1,\overline{A \setminus \mathcal{O}})\in \Lambda_j$.
The definition of $d_j$ and \eqref{y3}-\eqref{y4} then imply
\[d = d_{j} \leq \max_{u \in \eta(1, \overline{A \setminus \mathcal{O}})} J(u) \leq d - \varepsilon<d,\]
a contradiction. Thus $\operatorname{genus}(K^d)\ge i+1$. If $i\ge1$, then $\operatorname{genus}(K^d)\ge2$; a finite symmetric set not containing zero has genus one, so $K^d$ must be infinite.
\end{proof}

\begin{proposition}
  \label{p2}
  We have the following assertions:
  \begin{enumerate}[$(a)$]
    \item For every $n\in\mathbb{N}^+$, $c_n$ is a critical value of $J$,
    \item $b_n\geq c_n\geq \mu$,
    \item $c_n\rightarrow \infty$ as $n\rightarrow\infty$.
  \end{enumerate}
\end{proposition}

\begin{proof}
$(a)$ Since $J$ is even, satisfies the Palais--Smale condition by Lemma \ref{J1}, and
\[
 \max_{\partial D_n}J(\gamma_{0n})<0<c_n,
\]
the standard symmetric minimax theorem \cite[Chapter 9]{Rabinowitz1986}, applied to the homotopy-stable class $\Gamma_n$, shows that $c_n$ is a critical value of $J$.

$(b)$ For every $\gamma\in\Gamma_n$, oddness gives $\gamma(0)=0$, whereas $\|\gamma(\sigma)\|>\rho$ on $\partial D_n$ by Proposition \ref{p1}. Hence,
by continuity,
\[ \gamma(D_n) \cap \{u \in H_r^2(\mathbb{R}^N) : \|u\| = \rho\} \neq \emptyset \]
for all \( \gamma \in \Gamma_n \) and $\rho$ given in Proposition \ref{p1}. It follows from Proposition \ref{p1} and $I\ge J$ that
\[
 b_n\ge c_n\ge\mu>0.
\]

$(c)$ Taking $m=n$, $h=\gamma$, and $Y=\emptyset$ in the definition of $\Lambda_n$ gives $\gamma(D_n)\in\Lambda_n$ for every $\gamma\in\Gamma_n$.
Therefore
\[
d_n \leq \max_{u \in \gamma(D_n)} J(u) = \max_{\sigma \in D_n} J(\gamma(\sigma)).
\]
Taking the infimum over all $\gamma \in \Gamma_n$, we obtain
\[
d_n \leq c_n,
\]
for all \(n \in \mathbb{N}^+\).  Proposition \ref{p3} $(ii)$ and \eqref{mini1} imply that $\{d_n\}$ is
nondecreasing,
\[d_1 \leq d_2 \leq \cdots \leq d_n \leq d_{n+1} \leq \cdots.\]

We claim that $d_n\to\infty$. Otherwise, $d_n\uparrow \overline{d}<\infty$. If $d_n=\overline{d}$ for all large $j$, Lemma \ref{lem49} implies $\operatorname{genus}(K^{\overline{d}})=\infty$. But by Lemma \ref{J1}, $K^{\overline{d}}$ is compact so $\operatorname{genus}(K^{\overline{d}})<\infty$. Thus $\overline{d}>d_n$ for all $n\in\mathbb{N}^+$. The set
\[
 \mathcal K:=\{u\in H_r^2(\mathbb R^N):J|_{H_r^2}'(u)=0,
 \ d_{k+1}\le J(u)\le \overline{d}\}
\]
is compact by the Palais--Smale condition in Lemma \ref{J1}, symmetric, and does not contain zero. Hence $q:=\operatorname{genus}(\mathcal K)<\infty$, and there is a
symmetric neighborhood $\mathcal U$ of $\mathcal K$, with $0\notin\overline{\mathcal U}$, such that
$\operatorname{genus}(\overline{\mathcal U})=q$. Let $i=\max\{q,k+1\}$. The deformation theorem \cite[Theorem A.4]{Rabinowitz1986} with $c=\overline{d}$, $\overline{\varepsilon}=\overline{d}-d_s$, and $\mathcal O=\overline{\mathcal U}$ yields an $\varepsilon$ and $\eta$ as usual such that
\[\eta(1,J^{\overline{d}+\varepsilon}\setminus\mathcal{O})\subset J^{\overline{d}-\varepsilon}.\]
Choose $n\in \mathbb{N}^+$ such that $d_n>\overline{d}-\varepsilon$ and $A\in \Lambda_{n+i}$ such that
\[\max_{A}J\leq \overline{d}+\varepsilon.\]
The same argument as in the proof of Lemma \ref{lem49} shows that $\overline{A \setminus \mathcal{O}}\in \Lambda_n$, $\eta(1,\overline{A \setminus \mathcal{O}})\in \Lambda_n$ and provides an odd deformation that fixes all boundary sets $\gamma_{0m}(\partial D_m)$ and $\eta(1,u)=u$. Therefore, by the choice of $d_n$,
\[ d_{n} \leq \max_{u \in \eta(1, \overline{A \setminus \mathcal{O}})} J(u) \leq \overline{d} - \varepsilon,\]
a contradiction. Thus $d_n\to\infty$, and since $d_n\le c_n$, we conclude that $c_n\to\infty$.
\end{proof}

%Here, it seems difficult to show the Palais--Smale compactness condition for $I$ directly and it is a main difficulty in showing $b_n$ are critical values of $I$.
The lack of a direct Palais--Smale condition for $I$ is the main obstacle to proving that the levels $b_n$ are critical values.
To overcome this difficulty, introduce the augmented functional
%Here, we introduce an auxiliary functional
$\widetilde{I}\in C^1(\mathbb{R}\times H_r^2)$ by
 \begin{align*}\label{func2}
\widetilde{I}(s,v):=&I\bigl(v(e^{-s}\,\cdot)\bigr)\notag\\
=&\frac{1}{2}\left[e^{(N-4)s}\rn|\Delta v|^2dx+e^{(N-2)s}\rn |\nabla v|^2dx+e^{Ns}\rn v^2dx\right]\nonumber\\
  &+\frac{e^{(N-2)s}}{2}\rn v^2|\nabla v|^2dx-e^{Ns}\rn F(v)dx.
 \end{align*}

Now, we can define the minimax values $\widetilde{b}_n$ for $\widetilde{I}$ by
\begin{equation*}
\widetilde{b}_n = \inf_{\widetilde{\gamma} \in \widetilde{\Gamma}_n} \max_{\sigma \in D_n} \widetilde{I}(\widetilde{\gamma}(\sigma)),
\end{equation*}
where
\begin{equation*}
\widetilde{\Gamma}_n = \left\{ \widetilde{\gamma} \in C(D_n, \mathbb{R}\times H_r^2(\mathbb{R}^N)) \;\middle|\;
\begin{aligned}
&\widetilde{\gamma }=(s,\eta), s: D_n\rightarrow\mathbb{R} \text{\quad is even,}\\
&\eta: D_n\rightarrow H_r^2(\mathbb{R}^N)\text{\quad is odd,\quad}\widetilde{\gamma}|_{\partial D_n} =(0, \gamma_{0n})
\end{aligned}
\right\}.
\end{equation*}
%Observe that $\widetilde{\gamma}=(0,\gamma)\in\widetilde{\Gamma}_n$, for any $\gamma\in \Gamma_n$. On the other hand, for any $\widetilde{\gamma}=(s,\eta)\in\widetilde{\Gamma}_n $, $\gamma:=\eta(\cdot)(e^{-s(\cdot)}x)\in \Gamma_n$

\begin{lemma}
  \label{lem3}
  $\widetilde{b}_n=b_n$ for all $n\in\mathbb{N}^+$.
\end{lemma}

\begin{proof}
  If $\gamma\in\Gamma_n$, then the map $\sigma\mapsto(0,\gamma(\sigma))$ belongs to $\widetilde{\Gamma}_n$. Thus, by the definitions of $b_n$ and $\widetilde{b}_n$, and by $\widetilde{I}(0,u)=I(u)$, we obtain
  \begin{align}
    \label{s3}
    \max_{\sigma \in D_n} I(\gamma(\sigma)) = \max_{\sigma \in D_n}\widetilde{I} (0,\gamma(\sigma))\text{.}
  \end{align}
 This implies that $\widetilde{b}_n\leq b_n$. On the other hand, if $\widetilde{\gamma}(\sigma)=(s(\sigma),\eta(\sigma))\in \widetilde{\Gamma}_n$, we set
  \[\gamma(\sigma)=\eta(\sigma)(e^{-s(\sigma)}x)\text{.}\]
  Since $s$ is even and $\eta$ is odd, $\gamma(-\sigma)=-\gamma(\sigma)$; on $\partial D_n$, $s=0$ and $\eta=\gamma_{0n}$. Hence $\gamma\in\Gamma_n$, and $\widetilde{I}(\widetilde{\gamma}(\sigma))=I(\gamma(\sigma))$ for all $\sigma\in D_n$. Therefore,
  \[
\max_{\sigma \in D_n} \widetilde{I} (\widetilde{\gamma}(\sigma)) = \max_{\sigma \in D_n} I(\gamma(\sigma))
\]
and then
\[
\widetilde{b}_{n} = \inf_{\widetilde{\gamma} \in \widetilde{\Gamma}_n} \max_{\sigma \in D_n} \widetilde{I}(\widetilde{\gamma}(\sigma)) \geq \inf_{\gamma \in \Gamma_n} \max_{\sigma \in D_n} I(\gamma(\sigma)) = b_n.
\]
\end{proof}

We now prove that, for every $n\in\mathbb{N}^+$, $b_n$ is a critical value of $I$.

\begin{lemma}
  \label{lem4}
  For each fixed $n\geq 1$, there exists a sequence $\{v_j\}_{j\geq1}\subset H_r^2(\mathbb{R}^N)$ such that
  \begin{enumerate}[$(1)$]
    \item $I(v_j)\rightarrow b_n$;
    \item $\|I|_{H_r^2}'(v_j)\|_{(H_r^2)'}\rightarrow 0$;
    \item $P(v_j)\rightarrow 0$.
  \end{enumerate}
\end{lemma}

\begin{proof}
  By definition of $b_n$, for any $j\geq 1$ there exists $\gamma_j\in\Gamma_n$ such that
  \[\max_{\sigma \in D_n} I(\gamma_j(\sigma))\leq b_n+\frac{1}{j^2}\text{.}\]
  Setting $\widetilde{\gamma}_j = (0, \gamma_j) \in \widetilde{\Gamma}_n$, we deduce from \eqref{s3} and Lemma \ref{lem3} that
  \[\max_{\sigma \in D_n}\widetilde{I} (\widetilde{\gamma}_j(\sigma))\leq b_n+\frac{1}{j^2}\text{.}\]
  Since $\widetilde I(s,-u)=\widetilde I(s,u)$, the class $\widetilde\Gamma_n$ is stable under equivariant deformations associated with the involution $(s,u)\mapsto(s,-u)$. Applying the corresponding equivariant minimax principle (the equivariant form of \cite[Theorem 2.8]{MR1400007}), we obtain a sequence $\{(s_j, u_j)\} \subset \mathbb{R} \times H_r^2(\mathbb{R}^N)$ satisfying the following properties as $j \to \infty$,
    \begin{enumerate}[$(b_1)$]
      \item $\widetilde{I}(s_j,u_j)\rightarrow b_n$;
      \item  $\operatorname{dist}((s_j,u_j),\widetilde{\gamma}_j(D_n))\rightarrow 0$;
      \item $\|\widetilde{I}'(s_j,u_j)\|_{E_r'}\rightarrow 0$,
    \end{enumerate}
  where $E_r:=\mathbb{R}\times H_r^2(\mathbb{R}^N)$ and $E_r'$ denotes the dual space of $E_r$. Because  $\widetilde{\gamma}_j(D_n)\subset \{0\}\times H_r^2$, $(b_2)$ implies $s_j\rightarrow 0$. For any $(\mu,\omega)\in\mathbb{R}\times H_r^2(\mathbb{R}^N)$, we observe
\begin{eqnarray}\label{s21}
   (I\circ \varphi)'(s_j,u_j)[\mu,\omega]=I^\prime(\varphi(s_j,u_j))[\varphi(s_j,\omega)]+P(\varphi(s_j,u_j))\mu\rightarrow0.
\end{eqnarray}
Taking $\mu=1$ and $\omega=0$ in $(b_3)$, we get
\[P(\varphi(s_j,u_j))\rightarrow 0\text{,\quad as\quad}j\rightarrow\infty\text{.}\]
Set $w_j:=\varphi(s_j,u_j)$. Then, as $j\rightarrow\infty$,
\[I(w_j)\rightarrow b_n\text{,\quad and \quad}P(w_j)\rightarrow 0\text{.}\]
Renaming $w_j$ as $v_j$, we obtain $(1)$ and $(3)$. To prove $(2)$, let $\omega\in H_r^2(\mathbb R^N)$ and define
\[
       \widetilde\omega_j(x)=\omega(e^{s_j}x).
\]
Then $\varphi(s_j,\widetilde\omega_j)=\omega$. Since $s_j\to0$, the scaling operators
$\omega\mapsto\widetilde\omega_j$ are uniformly bounded on $H_r^2(\mathbb R^N)$. Therefore,
\[
 |\langle I|_{H_r^2}'(v_j),\omega\rangle|
 =|\widetilde I'(s_j,u_j)[(0,\widetilde\omega_j)]|
 \le o_j(1)\|\widetilde\omega_j\|
 \le o_j(1)\|\omega\|.
\]
Hence $\|I|_{H_r^2}'(v_j)\|_{(H_r^2)'}\to0$, and Lemma \ref{lem4} is proved.

\end{proof}

\begin{proof}[\textbf{Proof of Theorem \ref{th2}.}]
Fix $n\in\mathbb N^+$ and let $\{v_j\}$ be the sequence obtained in Lemma \ref{lem4}:
\[
  I(v_j)\to b_n>0,
  \qquad \|I|_{H_r^2}'(v_j)\|_{(H_r^2)'}\to0,
  \qquad P(v_j)\to0 .
\]
Repeating the estimates \eqref{t1}--\eqref{t2}, with $v_j$ in place of $u_n$, shows that
$\{v_j\}$ is bounded in $H_r^2(\mathbb R^N)$. Therefore, after passing to a subsequence,
\[
 v_j\rightharpoonup v\quad\hbox{in }H_r^2(\mathbb R^N),\qquad
 v_j\to v\quad\hbox{in }L^q(\mathbb R^N)\quad\hbox{for every }q\in(2,2_*),
\]
and $v_j(x)\to v(x)$ a.e. in $\mathbb R^N$. By \eqref{f2}, H\"older's inequality, the boundedness of $\{v_j\}$ in $L^2\cap L^{2_*}$, and $v_j\to v$ in $L^q$, $q\in(2,2_*)$, we obtain, after first taking $j\to\infty$ and then $\varepsilon\to0$,
\[
    \int_{\mathbb R^N}f(v_j)(v_j-v)\,dx\to0.
\]
Similar to the proof of \eqref{g111}, we also have
\[
    \int_{\mathbb R^N}f(v_j)v_j\,dx\to\int_{\mathbb R^N}f(v)v\,dx.
\]
Passing to the limit in $\langle I|_{H_r^2}'(v_j),\phi\rangle$ for radial $\phi\in C_c^\infty(\mathbb R^N)$ and then using density, we obtain $I|_{H_r^2}'(v)=0$. Consequently,
\begin{align*}
 \|v_j\|^2+2\int_{\mathbb R^N}v_j^2|\nabla v_j|^2dx
 =&\int_{\mathbb R^N}f(v_j)v_jdx+o(1)\\
 =&\int_{\mathbb R^N}f(v)vdx\\
 =&\|v\|^2+2\int_{\mathbb R^N}v^2|\nabla v|^2dx.
\end{align*}
As in the proof of Lemma \ref{J1}, weak lower semicontinuity of the quadratic part and \eqref{qslinear}, together with equality of the sums above, yields convergence of the quadratic norms and hence $v_j\to v$ strongly in $H_r^2(\mathbb R^N)$. Hence
\[
     I(v)=b_n,
     \qquad I|_{H_r^2}'(v)=0.
\]
Thus, by Palais principle of symmetric criticality, we obtain
%Since $I$ is invariant under the orthogonal action of $H_r^2(\mathbb R^N)$ is its fixed-point space, the principle of symmetric criticality gives 
$I'(v)=0$. Since $b_n\ge c_n\ge\mu>0$, this critical point $v$ is nontrivial. Finally, by $(b)$ and $(c)$ of Proposition \ref{p2}, we obtain $b_n\to\infty$. Therefore, equation \eqref{eq1} has infinitely many nontrivial radially symmetric solutions.
\end{proof}

\section*{Acknowledgments}
%\vspace{-0.5em}
L. F. Yin was supported by the Natural Science Foundation of Sichuan Province (No. 2024NSFSC1343) and the National Natural Science Foundation of China (No. 12401140). F. Wang was supported by the Fundamental Research Funds for the Central Universities (2682026ZTPY010), China.

\bigskip
\noindent{\bf Data availability}\\
No data were used for the research described in this article.

\bigskip

%\noindent{\bf Data Availability Statements}\\
%Data sharing not applicable to this article as no datasets were generated or analyzed during the current study.

\bigskip

\noindent{\bf Conflict of interest}\\
The authors declare that they have no competing interests.

%\bigskip
%
%\noindent{\bf Authors' contributions}\\
%The authors contributed equally to this work. All authors read and approved the final manuscript.

\bigskip

%\bibliographystyle{aml}
%\bibliography{autnoumous}

\begin{thebibliography}{10}

\bibitem{Alarcon2025}
\newblock S.~Alarc\'on, J.~Faya and C.~Rey,
\newblock Concentration on the boundary and sign-changing solutions for a
  slightly subcritical biharmonic problem,
\newblock \emph{J. Differential Equations}, \textbf{437} (2025), Paper No.
  113285, 49.
%\newblock \urlprefix\url{https://doi.org/10.1016/j.jde.2025.113285}.

\bibitem{MR695535}
\newblock H.~Berestycki and P.-L. Lions,
\newblock Nonlinear scalar field equations. {I}. {E}xistence of a ground state,
\newblock \emph{Arch. Rational Mech. Anal.}, \textbf{82} (1983), 313--345.
%\newblock \urlprefix\url{https://doi.org/10.1007/BF00250555}.

\bibitem{MR3276713}
\newblock S.~Chen, J.~Liu and X.~Wu,
\newblock Existence and multiplicity of nontrivial solutions for a class of
  modified nonlinear fourth-order elliptic equations on {$\mathbb{R}^N$},
\newblock \emph{Appl. Math. Comput.}, \textbf{248} (2014), 593--601.
%\newblock \urlprefix\url{https://doi.org/10.1016/j.amc.2014.10.021}.

\bibitem{MR3591225}
\newblock B.~Cheng and X.~Tang,
\newblock High energy solutions of modified quasilinear fourth-order elliptic
  equations with sign-changing potential,
\newblock \emph{Comput. Math. Appl.}, \textbf{73} (2017), 27--36.
%\newblock \urlprefix\url{https://doi.org/10.1016/j.camwa.2016.10.015}.



\bibitem{Fernandez2022}
\newblock A.~J. Fern\'andez, L.~Jeanjean, R.~Mandel and M.~Mari\c{s},
\newblock Non-homogeneous {G}agliardo-{N}irenberg inequalities in
  {$\mathbb{R}^N$} and application to a biharmonic nonlinear {S}chr\"odinger
  equation,
\newblock \emph{J. Differential Equations}, \textbf{330} (2022), 1--65.
%\newblock \urlprefix\url{https://doi.org/10.1016/j.jde.2022.04.037}.

\bibitem{MR330754}
\newblock J.~Frehse,
\newblock Zum {D}ifferenzierbarkeitsproblem bei {V}ariationsungleichungen
  h\"{o}herer {O}rdnung,
\newblock \emph{Abh. Math. Sem. Univ. Hamburg}, \textbf{36} (1971), 140--149.
%\newblock \urlprefix\url{https://doi.org/10.1007/BF02995917}.

\bibitem{Hirata2010}
\newblock J.~Hirata, N.~Ikoma and K.~Tanaka,
\newblock Nonlinear scalar field equations in {$\mathbb{R}^N$}: mountain pass and symmetric mountain pass approaches,
\newblock \emph{Topol. Methods Nonlinear Anal.}, \textbf{35} (2010), 253--276.

\bibitem{MR1430506}
\newblock L.~Jeanjean,
\newblock Existence of solutions with prescribed norm for semilinear elliptic
  equations,
\newblock \emph{Nonlinear Anal.}, \textbf{28} (1997), 1633--1659.
%\newblock \urlprefix\url{https://doi.org/10.1016/S0362-546X(96)00021-1}.

\bibitem{MR1084570}
\newblock A.~C. Lazer and P.~J. McKenna,
\newblock Large-amplitude periodic oscillations in suspension bridges: some new
  connections with nonlinear analysis,
\newblock \emph{SIAM Rev.}, \textbf{32} (1990), 537--578.
%\newblock \urlprefix\url{https://doi.org/10.1137/1032120}.

\bibitem{MR0259693}
\newblock J.-L. Lions,
\newblock \emph{Quelques m\'{e}thodes de r\'{e}solution des probl\`emes aux
  limites non lin\'{e}aires},
\newblock Dunod; Gauthier-Villars, Paris, 1969.

\bibitem{MR3348950}
\newblock H.~Liu and H.~Chen,
\newblock Least energy nodal solution for quasilinear biharmonic equations with
  critical exponent in {$\mathbb{R}^N$},
\newblock \emph{Appl. Math. Lett.}, \textbf{48} (2015), 85--90.
%\newblock \urlprefix\url{https://doi.org/10.1016/j.aml.2015.03.002}.

\bibitem{Liu2025}
\newblock J.~Liu, Z.~Zhang and Q.~Guan,
\newblock Multiplicity of normalized solutions to the biharmonic {S}chr\"odinger
  equation with mixed nonlinearities,
\newblock \emph{Complex Var. Elliptic Equ.}, \textbf{70} (2025), 1022--1047.
%\newblock \urlprefix\url{https://doi.org/10.1080/17476933.2024.2350991}.

\bibitem{MR3945610}
\newblock S.~Liu and Z.~Zhao,
\newblock Solutions for fourth-order elliptic equations on {$\mathbb{R}^N$}
  involving {$u\Delta(u^2)$} and sign-changing potentials,
\newblock \emph{J. Differential Equations}, \textbf{267} (2019), 1581--1599.
%\newblock \urlprefix\url{https://doi.org/10.1016/j.jde.2019.02.017}.

\bibitem{Lu2024}
\newblock Y.~Lu and X.~Zhang,
\newblock Existence of normalized positive solution of nonhomogeneous
  biharmonic {S}chr\"odinger equations: mass-supercritical case,
\newblock \emph{J. Fixed Point Theory Appl.}, \textbf{26} (2024), Paper No. 25,
  13.
%\newblock \urlprefix\url{https://doi.org/10.1007/s11784-024-01113-y}.

\bibitem{MR1050908}
\newblock P.~J. McKenna and W.~Walter,
\newblock Travelling waves in a suspension bridge,
\newblock \emph{SIAM J. Appl. Math.}, \textbf{50} (1990), 703--715.
%\newblock \urlprefix\url{https://doi.org/10.1137/0150041}.

\bibitem{Mederski2023}
\newblock J.~Mederski and J.~Siemianowski,
\newblock Biharmonic nonlinear scalar field equations,
\newblock \emph{Int. Math. Res. Not. IMRN}, (2023), 19963--19995.
%\newblock \urlprefix\url{https://doi.org/10.1093/imrn/rnac303}.

\bibitem{Rabinowitz1986}
\newblock P.~H. Rabinowitz,
\newblock \emph{Minimax Methods in Critical Point Theory with Applications to
  Differential Equations}, vol.~65 of CBMS Regional Conference Series in
  Mathematics,
\newblock Conference Board of the Mathematical Sciences, Washington, DC; by the
  American Mathematical Society, Providence, RI, 1986.
%\newblock \urlprefix\url{https://doi.org/10.1090/cbms/065}.

\bibitem{MR1400007}
\newblock M.~Willem,
\newblock \emph{Minimax Theorems}, vol.~24 of Progress in Nonlinear
  Differential Equations and Their Applications,
\newblock Birkh\"{a}user Boston, Inc., Boston, MA, 1996.
%\newblock \urlprefix\url{https://doi.org/10.1007/978-1-4612-4146-1}.

\bibitem{Yang2025}
\newblock M.~Yang, W.~Ye and X.~Zhang,
\newblock Nondegeneracy of positive solutions for a biharmonic {H}artree
  equation and its application,
\newblock \emph{J. Differential Equations}, \textbf{428} (2025), 796--849.
%\newblock \urlprefix\url{https://doi.org/10.1016/j.jde.2025.02.024}.

\bibitem{MR4291515}
\newblock L.-F. Yin and S.~Jiang,
\newblock Existence of nontrivial solutions for modified nonlinear fourth-order
  elliptic equations with indefinite potential,
\newblock \emph{J. Math. Anal. Appl.}, \textbf{505} (2022), Paper No. 125459.
%\newblock \urlprefix\url{https://doi.org/10.1016/j.jmaa.2021.125459}.

\bibitem{Zhang2025}
\newblock Z.~Zhang and Y.~Wang,
\newblock Normalized ground state solutions of the biharmonic {S}chr\"odinger
  equation with general mass-supercritical nonlinearities,
\newblock \emph{Appl. Math. Lett.}, \textbf{163} (2025), Paper No. 109415, 6.
%\newblock \urlprefix\url{https://doi.org/10.1016/j.aml.2024.109415}.

\end{thebibliography}
\providecommand{\href}[2]{#2}
\providecommand{\arxiv}[1]{\href{http://arxiv.org/abs/#1}{arXiv:#1}}
\providecommand{\url}[1]{\texttt{#1}}
\providecommand{\urlprefix}{URL }

\end{document}